 \input amstex

\documentstyle{amsppt} 
\magnification1200
\NoBlackBoxes
\pagewidth{6.5 true in}
\pageheight{9.25 true in}
\document

\topmatter
\title 
Extreme values of zeta and $L$-functions 
\endtitle
\author K. Soundararajan 
\endauthor 
\address{Department of Mathematics, 450 Serra Mall, Bldg. 380, Stanford University, 
Stanford, CA 94305-2125, USA}
  \endaddress
\email{ksound{\@}stanford.edu} \endemail
\thanks   The author is partially supported by the National Science Foundation (DMS 0500711) 
and the American Institute of Mathematics (AIM). 
\endthanks

\endtopmatter
\def\lam{\lambda}

\head 1. Introduction \endhead

\noindent In this paper we introduce a ``resonance" method to produce large
values of $|\zeta(\frac 12+it)|$ and large and small central values of $L$-functions. 

\proclaim{Theorem 1}  If $T$ is sufficiently large 
then there exists $t\in [T,2T]$ such that 
$$
|\zeta(\tfrac 12+it)| \ge \exp\Big( (1+o(1)) 
\frac{\sqrt{\log T}}{\sqrt{\log \log T}}\Big).
$$
Moreover uniformly in the range $3\le V\le \frac 15 \sqrt{\log T/\log \log T}$ 
we have that 
$$
\text{meas} \{t \in [T,2T] : \ \ |\zeta(\tfrac 12+it)| \ge e^V \} \gg \frac{T}{(\log T)^4}
 \exp\Big(-10\frac{V^2}{\log \frac{\log T}{8V^2\log V}}\Big).
$$
\endproclaim

The problem of obtaining large values of $|\zeta(\tfrac 12+it)|$ was 
first considered by E.C. Titchmarsh who showed that there exist arbitrarily large $t$ 
with $|\zeta(\tfrac 12+it)| \ge \exp( \log^{\alpha} t)$ for any $\alpha <\frac 12$ (see Theorem 8.12 of [15]).   
In [9] H.L. Montgomery proved that, assuming the Riemann 
Hypothesis, there exist arbitrarily large values $t$ such that 
$$
|\zeta(\tfrac 12+it)| \gg \exp\Big(\frac{1}{20} \frac{\sqrt{\log |t|}}{
\sqrt{\log \log |t|}}\Big).
$$
R. Balasubramanian and K. Ramachandra [2] proved a similar result 
unconditionally, showing that there are arbitrarily large $t$ such 
that 
$$
|\zeta(\tfrac 12+it)| \gg \exp\Big(B \frac{\sqrt{\log |t|}}{
\sqrt{\log \log |t|}}\Big),
$$
for some positive constant $B$.  Their method is based on obtaining lower bounds 
for the moments $\int_T^{2T} |\zeta(\frac 12+it)|^{2k} dt$.  Later Balasubramanian [1]
optimized their argument and found that $B=0.530\ldots$ is 
permissible.\footnote{The value of $B$ stated by him 
is $B=0.75\ldots$, but there appears to be a numerical error in 
the calculation.}   

As well as improving these results, our Theorem above 
suggests that there should be still larger values of $|\zeta(\tfrac 12+it)|$.   
A. Selberg (see [15]) has shown that as $t$ varies between $T$ and $2T$, $\log |\zeta(\tfrac 12+it)|$ 
has an approximately Gaussian distribution with mean $0$ and variance $\sim \tfrac 12\log \log T$.  
This suggests that the set of $t\in [T,2T]$ with $|\zeta(\tfrac 12+it)| \ge e^V$ should 
have measure about $T\exp(-V^2/\log \log T)$.  Our Theorem furnishes a lower bound for this 
measure of the type $T\exp(-cV^2/\log \log T)$ for some positive constant $c$ uniformly in 
the range $\log \log T \le V\le (\log T)^{\frac 12 -\delta}$ for any fixed $\delta>0$.   If this type of 
estimate were to persist for larger $V$, then we would expect to find values of $|\zeta(\tfrac 12+it)|$ 
of size $\exp(C\sqrt{\log T\log \log T})$ for some positive constant $C$.  
Indeed, recently D.W. Farmer, S.M. Gonek and 
C.P. Hughes [4] have suggested, based on several interesting heuristic considerations, 
that the maximum size of $|\zeta(\tfrac 12+it)|$ is about $\exp(C\sqrt{\log t\log \log t})$ 
with $C=1/\sqrt{2} +o(1)$.   

Complementing the lower bound of Theorem 1, we have 
shown in [14] that assuming the Riemann hypothesis 
$$
\text{meas} \{ t\in [T, 2T]: \ \ |\zeta(\tfrac 12+it)| \ge 
e^V\} \ll T \exp\Big(-(1+o(1))\frac{V^2}{\log \log T}  \Big), 
$$
in the range $10\sqrt{\log \log T} \le V =o(\log \log T \log \log \log T)$.   
When $V\ge \log \log T \log \log \log T$ this measure is  
$\ll T\exp(-cV\log V)$ for some positive constant $c$.   For a precise 
statement see the Theorem in [14].  

The main idea of our proof is to find a Dirichlet polynomial $R(t)= \sum_{n\le N} 
r(n) n^{-it}$ which `resonates' with $\zeta(\tfrac 12+it)$ and picks out its 
large values.  Precisely, we will compute the smoothed moments 
$$
M_1(R,T) = \int_{-\infty}^{\infty} |R(t)|^2 \Phi(\tfrac tT) dt, \ \ \ 
\text{and} \ \ \ 
M_2(R,T) = \int_{-\infty}^{\infty} \zeta(\tfrac 12+it) |R(t)|^2 \Phi(\tfrac tT)dt. \tag{1}
$$
Here $\Phi$ denotes a smooth, non-negative function, compactly supported in $[1,2]$, 
with $\Phi(y) \le  1$ for all $y$, and $\Phi(y)=1$ for $5/4\le y\le 7/4$.  
Plainly 
$$
\max_{T\le t\le 2T} |\zeta(\tfrac 12+it)| 
\ge \frac{|M_2(R,T)|}{M_1(R,T)}. 
$$
When $N\le T^{1-\epsilon}$ we may evaluate $M_1(R,T)$ and $M_2(R,T)$ 
easily.  These are two quadratic forms in the unknown coefficients 
$r(n)$, and the problem thus reduces to maximizing the ratio of 
these quadratic forms.  Solving this optimization problem we obtain 
Theorem 1.

This method generalizes readily to provide large and small 
central values in families of $L$-functions.  By contrast, 
the method of Montgomery does not appear to generalize to this 
situation.  Recently Z. Rudnick and the author ([11] and [12]) found 
a flexible method to obtain lower bounds for moments in many families 
of $L$-functions, but the bounds obtained here are superior.  

\proclaim{Theorem 2}   Let $X$ be large.  There 
exists a fundamental discriminant $d$ with 
$X\le |d| \le 2X$ such that 
$$
L(\tfrac 12,\chi_d) \ge \exp\Big( \Big(\frac{1}{\sqrt{5}}+o(1)\Big) \frac{\sqrt{\log X}}{\log \log X}\Big).
$$
 Moreover, there exists a fundamental discriminant $d$ with 
$X\le |d|\le 2X$ such that 
$$
|L(\tfrac12,\chi_d)|\le \exp\Big(- \Big(\frac{1}{\sqrt{5}}+o(1)\Big) \frac{\sqrt{\log X}}{\log \log X}\Big).
$$
 Here $\chi_d$ 
denotes the real primitive character associated to the fundamental discriminant $d$. 
\endproclaim

Previously, D.R. Heath-Brown (unpublished, see [6]) had shown that there 
 arbitrarily large fundamental discriminants $d$ 
such that 
$$
L(\tfrac 12,\chi_d) \gg \exp\Big(C\frac{\sqrt{\log |d|}}{\log \log |d|}\Big),
$$  
for some positive constant $C$.  Heath-Brown's idea was extended by J. Hoffstein 
and P. Lockhart [6] to prove a similar result for quadratic twists 
of any modular form.  Our method may be adapted to give an analogous improvement 
of their result.  

S.D. Chowla has conjectured that $L(\tfrac 12,\chi_d) >0$ for all fundamental discriminants $d$.
From [13] we know that $L(\tfrac 12,\chi_d) \neq 0$ for a large proportion ($\frac 78$) of fundamental 
discriminants $d$, and from [3] that $L(\tfrac 12,\chi_d) >0$ for a positive proportion of fundamental 
discriminants $d$.  Nevertheless, Theorem 2 tells us that there are very small values 
of $L(\tfrac 12,\chi_d)$, and arguing as in Theorem 1 we can also show that 
there are $\gg X\exp(-C\log X/\log \log X)$ discriminants $d$ with such a small 
value of $L(\tfrac 12,\chi_d)$.

We give one more example of this method.  Let $k$ denote an even integer
and let $H_k=H_k(1)$ denote the set of Hecke eigencuspforms
of weight $k$ for the full modular group $\Gamma=SL_2({\Bbb Z})$.
We write the Fourier expansion of $f \in H_k$ as
$$
f(z) = \sum_{n=1}^{\infty} \lam_f(n) n^{\frac{k-1}{2}} e(nz) 
$$
and normalize so that $\lam_f(1)=1$.  Note that, with 
our normalization, Deligne's bound reads $|\lam_f(n)| \le d(n)$ 
although we do not require it here.  Associated to $f$ is the $L$-function
$$
L(s,f) = \sum_{n=1}^{\infty } \frac{\lam_f(n)}{n^s}.  
$$
Recall that the sign of the functional equation for $L(s,f)$ is $i^k$.  When 
$k\equiv 2\pmod 4$ it follows that the central values $L(\tfrac 12,f)$ equal 
zero.   

\proclaim{Theorem 3}  For large $k\equiv 0 \pmod 4$ there 
exists $f\in H_k$ with 
$$
L(\tfrac 12, f) \ge \exp\Big( (1+o(1)) \frac{\sqrt{2\log k}}{\sqrt{\log 
\log k}}\Big).
$$
There also exists $f\in H_k$ with 
$$
L(\tfrac 12,f) \le \exp\Big(-(1+o(1)) \frac{\sqrt{2\log k}}
{\sqrt{\log \log k}}\Big).
$$ 
\endproclaim

In Theorem 1 we have attempted to optimize the large values of $|\zeta(\tfrac 12+it)|$ 
produced by our method.  In Theorems 2 and 3 we have tried instead to keep the 
exposition simple, and not pushed the method to its limit.  For example, with greater 
work we could take a longer resonator, allowing 
us  to replace the $1/\sqrt{5}$ appearing in Theorem 2 
with $\sqrt{3}$.   

The resonance method  is useful in producing omega results in other contexts 
as well.  For example, in work in progress A. Booker and the author have used it to obtain large character sums improving and simplifying the results in [5].   Using this method and adding 
their ideas, N. Ng [10] has obtained large and small values of $|\zeta^{\prime}(\rho)|$ 
where $\rho$ runs over zeros of $\zeta(s)$, and D. Milicevic [8] has obtained lower bounds 
for $L^{\infty}$ norms of eigenfunctions.   
 
{\sl Acknowledgments.} I am grateful to Greg Martin for a valuable suggestion.    

\head 2.  Large values of $|\zeta(\tfrac 12+it)|$: Proof of Theorem 1 \endhead

\noindent  Let $\Phi$ be a smooth function compactly supported in $[1,2]$, such that $0\le \Phi(t)\le 1$ 
always and $\Phi(t)=1$ for $t \in (5/4,7/4)$.  Let ${\hat \Phi}(y)=\int_{-\infty}^{\infty} \Phi(t)e^{-ity}dt$ 
denote 
the Fourier transform of $\Phi$.  Integrating by parts we note that
${\hat \Phi}(y) \ll_{\nu} |y|^{-\nu}$ for any integer $\nu\ge 1$.  

We first show how to evaluate the moments $M_1(R,T)$ and $M_2(R,T)$ 
defined in (1) when $N\le T^{1-\epsilon}$.  Observe that
$$
\align
\int_{-\infty}^{\infty} |R(t)|^2 \Phi(\tfrac tT) dt 
&= \sum_{m, n\le N} r(m)\overline{r(n)} \int_{-\infty}^{\infty} 
\Big(\frac{n}{m}\Big)^{it} \Phi\Big(\frac tT\Big) dt 
\\
&= T\sum_{m, n\le N} r(m)\overline{r(n)} {\hat \Phi}(T \log (m/n)).
\\
\endalign
$$
Since $N\le T^{1-\epsilon}$ we see that if $m\neq n$ then $T|\log (m/n)| \gg 
T^{\epsilon}$ so that ${\hat \Phi}(T\log (n/m)) \ll_{\epsilon} T^{-2}$ say.  
Therefore 
$$
\align
M_1(R,T) &= 
T{\hat \Phi}(0) \sum_{n\le N} |r(n)|^2 + O\Big( T^{-1} 
\Big(\sum_{n\le N} |r(n)|\Big)^2\Big) 
\\
&=T {\hat \Phi}(0) (1+O(T^{-1})) \sum_{n\le N} |r(n)|^2, \tag{2}
\\
\endalign
$$
by a simple application of Cauchy's inequality. 

Now consider 
$$
\int_{-\infty}^{\infty} |R(t)|^2 \sum_{k \le T} \frac{1}{k^{\frac 12+it}} 
\Phi\Big(\frac tT \Big) dt 
= T\sum_{m,n\le N} \sum_{k\le T} 
\frac{r(m)\overline{r(n)}}{\sqrt{k}} {\hat \Phi}(T \log (mk/n)).
$$
If $N\le T^{1-\epsilon}$ then for off-diagonal terms $mk\neq n$ we 
have ${\hat \Phi}(T\log (mk/n))\ll_\epsilon T^{-2}$.  
Thus the above equals
$$
\align
&{\hat \Phi}(0) T\sum_{mk=n \le N} \frac{r(m)\overline{r(n)}}{\sqrt{k}} 
+ O\Big(T^{-1} \sum_{k\le T} \frac{1}{\sqrt{k}} \Big(\sum_{n\le N} |r(n)|\Big)^2 
\Big)\\
 = &T{\hat \Phi}(0) \sum_{mk=n \le N} \frac{r(m)\overline{r(n)}}{\sqrt{k}} 
+ O\Big(T^{\frac 12} \sum_{n\le N} |r(n)|^2\Big). 
\\
\endalign
$$
Since $\zeta(\frac 12+it) = \sum_{k\le T} k^{-\frac 12-it} +O(T^{-\frac 12})$ 
for $T\le t\le 2T$ (see Theorem 4.11 of  [15]) we deduce that 
$$
M_2(R,T) 
= T {\hat \Phi}(0) \sum_{mk=n \le N} \frac{r(m)\overline{r(n)}}{\sqrt{k}} 
+ O\Big(T^{\frac 12} \sum_{n\le N} |r(n)|^2\Big).  \tag{3} 
$$

From (2) and (3) we glean that, if $N\le T^{1-\epsilon}$ then  
$$
\max_{T\le t\le 2T} |\zeta(\tfrac 12+it)| 
\ge (1+O(T^{-1})) \Big|\sum_{mk \le N} \frac{r(m)\overline{r(mk)}}{\sqrt k} 
\Big| \Big/ \Big(\sum_{n\le N }|r(n)|^2\Big) 
+O(T^{-\frac 12}). \tag{4}
$$ 
It remains to choose the resonator coefficients $r(n)$ so as to maximize this 
ratio.

\proclaim{Theorem 2.1}  For large $N$ we have 
$$
\max_{r} \Big|\sum_{mk \le N} \frac{r(m)\overline{r(mk)}}{\sqrt k} 
\Big| \Big/ \Big(\sum_{n\le N }|r(n)|^2\Big) 
=\exp\Big( \frac{\sqrt{\log N}}{\sqrt{\log \log N}} +O\Big(
\frac{\sqrt{\log N}}{\log \log N}\Big)\Big). 
$$
\endproclaim 

\demo{Proof of the lower bound of Theorem 2.1}  
We take $r(n)$ to be $f(n)$ where $f$ is a multiplicative function 
such that $f(p^k)=0$ for $k\ge 2$.  Let $L:=\sqrt{\log N\log \log N}$, 
and define $f(p)=L/(\sqrt{p}\log p)$ if $L^2 \le p\le \exp((\log L)^2)$, and 
$f(p)=0$ for all other primes $p$.  
Note that the denominator in our ratio is 
$$
\sum_{n\le N} f(n)^2 \le \sum_{n=1}^{\infty} f(n)^2 = 
\prod_{p} (1+f(p)^2). \tag{5}
$$

Now we need a lower bound for the numerator of our ratio.  
Below we make use of the observation that if $a_n$ is a
sequence of non-negative real numbers then for any $\alpha >0$ 
we have 
$$
\sum_{n>x} a_n \le x^{-\alpha} \sum_{n>x} a_n n^{\alpha} \le 
x^{-\alpha} \sum_{n=1}^{\infty} a_n n^{\alpha}. 
$$
This observation is often called `Rankin's trick.'
Thus the numerator of our ratio is, for any $\alpha >0$,
$$
\sum_{k\le N} \frac{f(k)}{\sqrt{k}} \sum\Sb n\le N/k\\ (n,k)=1\endSb f(n)^2 
=\sum_{k\le N} \frac{f(k)}{\sqrt{k}} \Big( 
\prod_{p \nmid k} (1+f(p)^2) + O\Big( \Big(\frac{k}{N}\Big)^{\alpha} 
\prod_{p\nmid k} (1+p^\alpha f(p)^2)\Big)\Big).
$$  
The error term above is plainly 
$$
O\Big( \frac{1}{N^{\alpha}} \prod_p (1+p^\alpha f(p)^2 
+ f(p)p^{\alpha -\frac 12})\Big), 
$$
while the main term is 
$$
\sum_{k\le N}  \frac{f(k)}{\sqrt{k}} \prod_{p \nmid k} (1+f(p)^2) 
= \prod_{p} (1+f(p)^2 + f(p)/\sqrt{p}) + O
\Big(\frac{1}{N^{\alpha}} \prod_{p} (1+f(p)^2 +f(p)p^\alpha/\sqrt{p})\Big).
$$
Thus our numerator is 
$$
\prod_{p} (1+f(p)^2 + f(p)/\sqrt{p}) + 
O\Big(\frac{1}{N^{\alpha}} \prod_p (1+p^\alpha f(p)^2 
+ f(p)p^{\alpha -\frac 12})\Big). \tag{6}
$$
Taking $\alpha =1/(\log L)^3$ we may see that 
the ratio of the error term in (6) to the 
main term there is 
$$
\ll \exp\Big(-\alpha \log N +\sum_{L^2 \le p\le \exp(\log^2 L)} 
(p^\alpha-1) \Big(\frac{L}{p\log p}+\frac{L^2}{p\log^2 p}\Big) 
\Big) \ll \exp\Big(-\alpha \frac{\log N}{\log \log N}\Big),
$$
with a little calculation using the prime number theorem.  
Thus for large $N$ the numerator of our ratio is at least 
$$
\frac 12 \prod_{p}\Big( 1+f(p)^2+ \frac{f(p)}{\sqrt{p}}\Big), \tag{7}
$$ 
and the lower bound of the Theorem follows from (5). 
\enddemo 

\demo{Proof of the upper bound of Theorem 2.1}  Define the 
multiplicative function $g$ by setting 
$g(p^k)=\min (1, L/(p^{k/2} \log p))$ where 
$L=\sqrt{\log N\log \log N}$ as above.  
Since $2|r(mk)r(m)| \le |r(mk)|^2/g(k) + g(k) |r(m)|^2$ we obtain that 
the numerator of our ratio is 
$$
\le \frac{1}{2} 
\sum_{km\le N} \frac{1}{\sqrt{k}} \Big(\frac{|r(mk)|^2}{g(k)} + 
g(k)|r(m)|^2\Big) 
=\frac12 
\sum_{n\le N} |r(n)|^2 \Big(\sum_{k\le N/n} \frac{g(k)}{\sqrt{k}} 
+ \sum_{k|n} \frac{1}{\sqrt{k}g(k)}\Big),
$$
with a little regrouping.  Note that 
$$
\align
\sum_{k\le N/n} \frac{g(k)}{\sqrt{k}} 
&\le \prod_{p} \Big( 1 +\frac{g(p)}{\sqrt{p} -1}\Big) \\
&\ll \exp\Big( \sum_{p\le \log N/\log \log N} \frac{1}{\sqrt{p}-1} + 
\sum_{p> \log N/\log \log N} \frac{L}{\sqrt{p}(\sqrt{p}-1)\log p} \Big) \\
&= \exp\Big( \frac{\sqrt{\log N}}{\sqrt{\log \log N}}+ O\Big( 
\frac{\sqrt{\log N}\log \log \log N}{(\log \log N)^{\frac 32}}\Big)\Big).\\
\endalign
$$
Further observe that for $n\le N$ 
$$
\sum_{k|n} \frac{1}{\sqrt{k}g(k)} \le 
\prod_{p^a \parallel n} \Big(1+\frac{a\log p}{L}\Big) \prod_{p|n} \Big( 1+ \frac{1}{\sqrt{p}-1}\Big).
$$
The first factor above is $\le \exp(\sum_{p^a \parallel n} (a\log p)/L) = n^{1/L} \le N^{1/L}$.  
The second factor is $\ll \exp(O(\sqrt{\log N}/\log \log N))$ by a simple calculation using 
the prime number theorem.  
The upper bound implicit in the Theorem follows. 
 
\enddemo

\demo{Proof of Theorem 1}  Using Theorem 2.1 in (4), and choosing $N=T^{1-\epsilon}$ we obtain immediately the first assertion of Theorem 1.   It remains now to establish the 
lower bound on the frequency with which large values are attained.  We have 
$$
T\log T \sim \int_T^{2T} |\zeta(\tfrac 12+it)|^2 dt 
\le \frac 12 T\log T + \int\Sb t \in [T,2T] \\ |\zeta(\frac 12+it)| \ge \sqrt{\frac 12\log T}\endSb 
|\zeta(\tfrac 12+it)|^2 dt, 
$$
so that 
$$
\align
T\log T & \ll \int \Sb t \in [T,2T] \\ |\zeta(\frac 12+it)| \ge \sqrt{\frac 12\log T}\endSb 
|\zeta(\tfrac 12+it)|^2 dt\\
& \ll \Big( \text{meas} \{t\in [T,2T]: |\zeta(\tfrac 12+it)| \ge 
\sqrt{\tfrac 12 \log T}\} \Big)^{\frac 12} \Big(\int_{T}^{2T} |\zeta(\tfrac 12+it)|^4 dt\Big)^{\frac 12}.\\
\endalign
$$
Since $\int_T^{2T} |\zeta(\tfrac 12+it)|^4 dt \asymp T(\log T)^4$ we conclude 
that
$$
 \text{meas} \{t\in [T,2T]: |\zeta(\tfrac 12+it)| \ge 
\sqrt{\tfrac 12 \log T}\}  \gg \frac{T}{(\log T)^2}, 
$$
which gives our desired lower bound when $3\le V\le \frac 12 \log(\frac 12\log T)$. 

For larger values of $V$, we use the resonator method with $N=T^{\frac 12-\epsilon}$.  
If  $2e^V M_1(R,T) \le |M_2(R,T)|$ (with the resonator $R$ still to be chosen) then 
$$
\align
|M_2(R,T) | 
&\le e^V M_1(R,T) +  \int_{\{t: |\zeta(\frac 12+it)| \ge e^V\} } |\zeta(\tfrac 12+it)| |R(t)|^2 
\Phi(\tfrac tT) dt.\\
\endalign
$$
Using Cauchy's inequality twice, we see that the integral above is  
$$
\le \Big(\text{meas} \{t \in [T,2T]:  \ |\zeta(\tfrac 12+it)| \ge e^V \} 
\Big)^{\frac 14} \Big( \int_T^{2T} |\zeta(\tfrac 12+it)|^4 dt \Big)^{\frac 14} 
\Big( \int_{-\infty}^{\infty} |R(t)|^4 \Phi(\tfrac tT) dt\Big)^{\frac 12}.  
$$
Therefore
$$
\text{meas} \{t \in [T,2T]:  \ |\zeta(\tfrac 12+it)| \ge e^V \}  
\gg\frac{ |M_2(R,T)|^4}{T \log^4 T} \Big( \int_{-\infty}^{\infty} |R(t)|^4 \Phi(\tfrac tT)dt \Big)^{-2}.
\tag{8}
$$

Let $A$ be large with $10A^2 \log A \le \log N$.  We choose the resonator coefficients 
$r(n)$ to be multiplicative, with $r(p^k)=0$ for $k\ge 2$ and 
$$
r(p)= \cases 
A/\sqrt{p} &\text{if } A^2 \le p \le N^{\frac{1}{2A^2}}\\ 
0&\text{otherwise}.\\
\endcases
$$ 
We use Rankin's trick and argue as in the proof of the lower bound of Theorem 2.1 
(taking now $\alpha=A^2/\log N$).   That gives
$$
\frac{|M_2(R,T)|}{M_1(R,T) }\ge 
\frac 12\prod_p \Big( 1+r(p)^2+ \frac{r(p)}{\sqrt{p}}\Big)(1+r(p)^2)^{-1} 
= \exp\Big( (1+o(1)) A \log \frac{\log N}{4A^2 \log A}\Big).  \tag{9}
$$
Further, 
$$
\int_{-\infty}^{\infty} |R(t)|^4 \Phi(\tfrac tT) dt 
= \sum_{a,b,c,d \le N} r(a)r(b)r(c)r(d)   T{\hat \Phi}(T\log \tfrac{ab}{cd}).
$$
If $N\le T^{\frac 12-\epsilon}$ then if $ab\neq cd$ then $|\log \frac{ab}{cd}| 
\gg T^{-1+\epsilon}$ so that ${\hat \Phi}(T\log \frac {ab}{cd}) \ll_{\epsilon} T^{-4}$, 
say.  
Since $r(n) \le 1$ for all $n$ we conclude that the off-diagonal terms $ab\neq cd$ 
contribute an amount  $\ll T^{-3} N^4 \ll T^{-1}$.  Thus 
$$
\align
\int_{-\infty}^{\infty} |R(t)|^4 \Phi(\tfrac tT) dt 
&= T{\hat \Phi}(0) \sum\Sb a,b,c,d \le N\\ ab=cd \endSb r(a)r(b)r(c)r(d)  
+O(T^{-1})\\
&\ll T \prod_p (1+ 4 r(p)^2 + r(p)^4). 
\\
\endalign
$$
Since $M_2(R,T) \gg T$, we conclude from (8) that  
$$
\align
\text{meas} \{t \in [T,2T]:  \ |\zeta(\tfrac 12+it)| \ge e^V \}  
&\gg \frac{T}{\log^4 T} \exp\Big(-4\sum_p r(p)^2 \Big) 
\\
&\gg \frac{T}{\log^4 T} \exp\Big( - 5 A^2 \log \frac{\log N}{4A^2 \log A}\Big). \\
\endalign
$$
We may choose $A \sim V (\log \frac{\log N}{4V^2\log V})^{-1}$ such that 
the RHS of (9) exceeds $2e^V$, and then the above estimate yields the 
bound claimed in Theorem 1. 

\enddemo

\head 3. Extreme values of quadratic Dirichlet $L$-functions: Proof of Theorem 2\endhead
\def\L{\fracwithdelims()}

\noindent For convenience, we restrict ourselves to fundamental discriminants 
of the form $8d$ where $d$ is an odd, squarefree number with $X/16 \le d\le X/8$.  
As before, we will consider the two moments 
$$
M_1(R,X) = \sum_{ X/16 \le d\le X/8 } \mu(2d)^2 R(8d)^2, \ \ \ 
M_2(R,X) = \sum_{X/16 \le d\le X/8} \mu(2d)^2 L(\tfrac 12,\chi_{8d}) R(8d)^2, 
$$
where 
$$
R(8d) = \sum_{n\le N} r(n) \L{8d}{n}, 
$$
is a resonator, whose coefficients $r(n)$ are real numbers to be chosen presently.  

\proclaim{Lemma 3.1}  The quantity $M_1(R,X)$ equals
$$
 \frac{X}{16\zeta(2)} 
\sum\Sb n_1,n_2\le N \\ n_1 n_2= \text{ odd 
square}\endSb r(n_1)r(n_2) \prod_{p|2n_1n_2} \L{p}{p+1} 
+ O\Big(X^{\frac 12+\epsilon} N^{\frac 12} \Big(\sum_{n\le N} |r(n)|\Big)^2\Big).
$$
\endproclaim
\demo{Proof} 
Expanding $R(8d)^2$ we see that 
$$
M_1(R,X) = \sum_{n_1, n_2 \le N} r(n_1) r(n_2) 
\sum_{X/16 \le d\le X/8} \mu(2d)^2 \L{8d}{n_1n_2}. \tag{10} 
$$
Let $n$ be an odd number and $z\ge 3$.  We record the following character 
sum estimate which may be obtained easily from the P{\' o}lya-Vinogradov inequality
(or see Lemma 3.1 of [12] for details).
If $n$ is not a perfect square then 
$$
\sum_{d\le z} \mu(2d)^2 \L{8d}{n} \ll z^{\frac 12} n^{\frac 14} \log (2n), \tag{11a}
$$
while if $n$ is a perfect square then 
$$
\sum_{d\le z} \mu(2d)^2 \L{8d}{n} = \frac{z}{\zeta(2)} \prod_{p|2n} \L{p}{p+1} + O (z^{\frac 12+\epsilon} 
n^{\epsilon}). \tag{11b}
$$
The Lemma follows upon using (11a,b) in (10).
\enddemo


To evaluate $M_2(R,X)$ we will use (11a,b) along with a standard ``approximate 
functional equation." The approximate functional equation we need states that for an 
odd, positive, square-free number $d$ we have
$$
L(\tfrac 12,\chi_{8d}) = 2\sum_{n=1}^{\infty} \frac{\chi_{8d}(n)}{\sqrt n}W\L{n\sqrt{\pi}}{\sqrt{8d}}
$$
where the weight $W$ is defined by 
$$
W(\xi) =\frac{1}{2\pi i} \int_{(c)} \frac{\Gamma(\frac s2+\frac 14)}{\Gamma(\frac 14)} 
\xi^{-s} \frac{ds}{s},
$$
and the integral is over a vertical line $c-i\infty$ to $c+i\infty$ with $c>0$.  The weight $W(\xi)$ 
is smooth and satisfies $W(\xi) =1 +O(\xi^{\frac 12-\epsilon})$ for small $\xi$, and $W(\xi) \ll e^{-\xi}$ 
for large $\xi$.  Moreover the derivative $W^{\prime}(\xi)$ satisfies $W^{\prime}(\xi) \ll \xi^{\frac 12-\epsilon} e^{-\xi}$.  These facts are easily established; for details see Lemmas 2.1 and 2.2 of [13], 
or Lemma 3.2 of [12].  

\proclaim{Lemma 3.2}  The quantity $M_2(R,X)$ equals 
$$
\align
\frac{X}{8\zeta(2)} \sum\Sb n_1, n_2\le N \endSb 
&r(n_1)r(n_2) \sum\Sb n \\ nn_1n_2 = \text{ odd square} \endSb 
\frac{1}{\sqrt{n}} \prod_{p|2nn_1n_2} \L{p}{p+1} 
\int_1^2 W\L{n\sqrt{2\pi}}{\sqrt{Xt}} dt  
\\
&+ O\Big( X^{\frac 78+\epsilon} N^{\frac 12} \Big(\sum_{n\le N} |r(n)| \Big)^2 \Big).
  \\
\endalign
$$
\endproclaim 
\demo{Proof}  Expanding $R(8d)^2$, and using the approximate functional 
equation, we have that
$$
M_2(R,X) = 2\sum_{n_1, n_2\le N} 
r(n_1) r(n_2) \sum_{n=1}^{\infty} \frac{1}{\sqrt{n}} 
\sum_{X/16\le d\le X/8}\mu(2d)^2 \L{8d}{nn_1n_2} W\L{n\sqrt{\pi}}{\sqrt{8d}}.
$$
By (11a,b) and partial summation we see that if $nn_1n_2$ is not an 
odd square then 
$$
\sum_{X/16\le d\le X/8}\mu(2d)^2 \L{8d}{nn_1n_2} W\L{n\sqrt{\pi}}{\sqrt{8d}} \ll X^{\frac 12} (nn_1n_2)^{\frac 14+\epsilon} e^{-n/\sqrt{X}},
$$
while if $nn_1n_2$ is an odd square that sum over $d$ is 
$$
\frac{X}{16\zeta(2) }\prod_{p|2nn_1n_2} \L{p}{p+1} \int_1^2 W\L{n\sqrt{2\pi}}{\sqrt{Xt}} dt 
+ O(X^{\frac 12+\epsilon} e^{-n/\sqrt{X}} ).
$$
The errors above contribute to $M_2(R,X)$ an amount
$$
\ll X^{\frac 12+\epsilon} N^{\frac 12} \Big( \sum_{\ell \le N} |r(\ell)|\Big)^2  \sum_{n=1}^{\infty} 
\frac{n^{\frac 14+\epsilon}}{\sqrt{n}} e^{-n/\sqrt{X}} 
\ll X^{\frac 78 +\epsilon} N^{\frac 12} \Big(\sum_{n\le N} |r(n)|\Big)^2.
$$ 
The Lemma follows.
\enddemo



\proclaim{Proposition 3.3}  Let $N\le X^{\frac 1{20}-\epsilon}$ be large.  Set $L=\sqrt{\log N\log \log N}$ 
and choose the resonator coefficients $r(n)$ to be $\mu(n)f(n)$ where $f$ is a multiplicative 
function with $f(p)=L/(\sqrt{p}\log p)$ for $L^2 \le p \le \exp((\log L)^2)$ and 
$f(p)=0$ for all other primes.\footnote{Thus $f$ is the function appearing in the 
proof of the lower bound in Theorem 2.1.}  Then 
$$
M_1(R,X) \sim \frac{X}{24\zeta(2)}\prod_p (1+f(p)^2),
$$
and 
$$
M_2(R,X) \sim  C_1X(\log X) \prod_p \Big( 1+f(p)^2 -2\frac{f(p)}{\sqrt{p}}\Big), 
$$
where $C_1$ is an absolute positive constant.   
\endproclaim

Taking $N=X^{\frac 1{20}-\epsilon}$, a little calculation shows that 
$$
\frac{M_2(R,X)}{M_1(R,X)} =\exp\Big( -(2+o(1)) \frac{\sqrt{\log N}}{\sqrt{\log \log N}}\Big) 
= \exp\Big( -\Big(\frac{1}{\sqrt{5}}+o(1)\Big)\frac{\sqrt{\log X}}{\sqrt{\log \log X}}\Big).
$$
 This demonstrates the existence of the small values claimed in Theorem 2.   
To find large values we take $r(n)=f(n)$ in Proposition 3.3, and argue in an identical 
manner.

\demo{Proof of Proposition 3.3} 
Lemma 3.1 gives 
$$
M_1(R,X) = \frac{X}{16\zeta(2)} \sum_{n\le N} \mu(n)^2 f(n)^2 \prod_{p|2n} \L{p}{p+1} 
+ O(X^{\frac 12+\epsilon} N^{\frac 52}).
$$
Rankin's trick shows that for any $\alpha >0$ 
$$
\sum_{n\le N} \mu(n)^2 f(n)^2 \prod_{p|2n} \L{p}{p+1}  
= \frac 23 \prod_p \Big(1+ f(p)^2 \frac{p}{p+1} \Big) + 
O\Big( N^{-\alpha} \prod_p \Big(1 +f(p)^2 p^{\alpha} \frac{p}{p+1} \Big) \Big).
$$
Choosing $\alpha= 1/(\log L)^3$ (as in Theorem 2.1) we find that the 
ratio of the error term above to the main term is $\ll \exp(-\alpha \log N/\log \log N)$. 
When $N \le X^{\frac 15-\epsilon}$ we conclude that 
$$
M_1(R,X) \sim \frac{X}{24\zeta(2)} \prod_{p} \Big(1 +f(p)^2 \frac{p}{p+1}\Big)
\sim \frac{X}{24\zeta(2)} \prod_{p} (1+f(p)^2).
$$
This proves the first assertion of the Proposition. 

Now we turn to $M_2(R,X)$.   We use Lemma 3.2, and note that when 
$N \le X^{\frac 1{20}-\epsilon}$  the remainder term there is $O(X^{1-\epsilon})$.   
Consider the main term in the asymptotic formula of Lemma 3.2.  To 
analyze this we write $n_1=ar$ and $n_2=as$ where $a=(n_1,n_2)$ so that $(r,s)=1$; from our 
choice of the coefficients $r(n)$, we also have that $(a,r)=(a,s)=1$.  With this notation, we may write 
the variable $n$ in Lemma 3.2 as $rsm^2$ for some odd integer $m$.  Thus the main term in Lemma 
3.2 equals 
$$
\align
\frac{X}{12 \zeta(2)} \sum\Sb a, r, s \\ ar, as\le N\\ (a,r)=(a,s)=(r,s)=1\endSb 
&\mu(a)^2 f(a)^2 \frac{\mu(r)f(r) \mu(s)f(s) }{\sqrt{rs}} 
\\
&\times \sum\Sb m\text{ odd} \endSb \frac{1}{m} \prod_{p|arsm} \L{p}{p+1} 
\int_1^2 W\Big( \frac{rsm^2 \sqrt{2\pi}}{\sqrt{Xt}}\Big) dt. 
\tag{12} \\
\endalign
$$
We now evaluate the sum over $m$ above.  Recalling the definition of $W(\xi)$ we 
may express that sum as 
$$
\frac{1}{2\pi i} \int_{(c)} \frac{\Gamma(\frac w2 +\frac 14)}{\Gamma(\frac 14)} \Big( \frac {X}{rs\sqrt{2\pi}}\Big)^{\frac{w}{2}} 
\Big(\int_1^2 t^{\frac w2} dt \Big) \sum_{m \text{ odd} } \frac{1}{m^{1+w}} \prod_{p|arsm} \L{p}{p+1} 
\frac{dw}{w}, \tag{13}
$$
where the integral is over the line from $c-i\infty$ to $c+i\infty$ with $c>0$.  A little calculation 
allows us to write the sum over $m$ above as 
$$
\zeta(1+w) (1-2^{-(1+w)}) \prod_{p|ars} \L{p}{p+1} 
\prod_{p\nmid 2 ars} \Big( 1-\frac{1}{p^{1+w}(p+1)}\Big). 
$$
We insert this in (13) and move the line of integration to $\text{Re } w= -\frac{1}{2}+\epsilon$.  
In view of the rapid decay of $\Gamma(\frac w2 +\frac 14)$, the integral on that 
line is $\ll (X/rs)^{-\frac 14+\epsilon}$, and therefore (13) equals 
$$
\align
\prod_{p|ars} \L{p}{p+1} {{\mathop{\text{Res}}}\atop{w=0}} \frac{\Gamma(\tfrac w2+\frac 14)}{\Gamma(\frac 14)} 
\Big(\frac{X}{rs\sqrt{2\pi}}\Big)^{\frac {w}{2}}& \Big(\int_{1}^{2} t^{\frac w2}dt\Big) 
\frac{\zeta(1+w)}{w}  (1-2^{-1-w})  \\
&\times \prod_{p\nmid 2ars} \Big(1-\frac{1}{p^{1+w}(p+1)} \Big) + O(X^{-\frac 14+\epsilon} (rs)^{\frac 14}).
\\
\endalign
$$
Computing the residue, we see that the above equals 
$$
\frac{1}{2} \prod_{p|ars}\L{p}{p+1} 
\prod_{p\nmid 2ars} \Big( 1-\frac{1}{p(p+1)}\Big) 
\Big( \log \frac{X}{rs} + C - \sum_{p|ars} \frac{\log p}{p(p+1)} \Big) + O(X^{-\frac 14+\epsilon}(rs)^{\frac 14}),
$$ 
for a suitable absolute constant $C$.   Using this in (12) we conclude that 
for $N\le X^{\frac 1{20}-\epsilon}$ 
$$
\align
M_2(R,X) = C_1 X &\sum\Sb a, r, s \\ ar, as\le N\\ (a,r)=(a,s)=(r,s)=1\endSb 
\mu(a)^2 f(a)^2 h(a) \frac{\mu(r)f(r)h(r)}{\sqrt{r}} \frac{\mu(s)f(s)h(s)}{\sqrt{s}} \\
&\times 
\Big(\log \frac{X}{rs} +C -\sum_{p|ars} \frac{\log p}{p(p+1)} \Big) + O(X^{1-\epsilon}),
\tag{14} \\
\endalign
$$
where $C_1$ is an absolute positive constant, and $h$ is a completely 
multiplicative function defined by $h(p)= p^2/(p^2+p-1)$. 

To simplify (14) further, we first extend the summations over $a$, $r$, and $s$ 
to run over all integers, and then use Rankin's trick to estimate the tails.   The extended sum equals 
$$
\align
\sum\Sb a, r, s \\ (a,r)=(a,s)=(r,s)=1\endSb 
\mu(a)^2 f(a)^2h(a)& \frac{\mu(r)f(r)h(r)}{\sqrt{r}} \frac{\mu(s)f(s)h(s)}{\sqrt{s}} 
\\
&\times \Big( \log X+ C- \sum_{p|rs} \log p -\sum_{p|ars} \frac{\log p}{p(p+1)}\Big).
\\
\endalign
$$
By multiplicativity this is seen to be
$$
\align
\prod_p &\Big( 1+ f(p)^2 h(p)- 2\frac{f(p)h(p)}{\sqrt{p}}\Big) 
\\
&\times \Big( \log X +C + \sum\Sb \ell \text{ prime}\endSb 
\log \ell \frac{2f(\ell)h(\ell)(1+1/(\ell(\ell+1)))/\sqrt{\ell} -f(\ell)^2 h(\ell)/(\ell(\ell+1))}{1+f(\ell)^2h(\ell)
-2f(\ell)h(\ell)/\sqrt{\ell}}\Big).\\
\endalign
$$
Since $h(p)=1+O(1/p)$ we may further simplify the above to 
$$
\sim (\log X)\prod_{p} \Big( 1+f(p)^2 -2\frac{f(p)}{\sqrt{p}}\Big). \tag{15}
$$
It remains to bound the error incurred upon extending the sums to infinity.  By symmetry 
we may suppose that $ar >N$, and we wish to estimate 
$$
\sum\Sb a,r,s\\ ar>N\endSb f(a)^2 \frac{f(r)f(s)}{\sqrt{rs}} (\log X+ \log r) 
\ll \prod_{p} \Big(1+\frac{f(p)}{\sqrt{p}}\Big) 
\sum\Sb ar>N \endSb f(a)^2 \frac{f(r)}{\sqrt{r}} (\log X+ \log r).
$$
As before we will use Rankin's trick with $\alpha= 1/(\log L)^3$.  If $ar>N$ then 
we have $(\log X+\log r) \ll (\log X) N^{-\alpha}(ar)^{\alpha}$.   Therefore, our desired 
quantity is 
$$
\align
&\ll (\log X) \prod_p \Big(1+ \frac{f(p)}{\sqrt{p}}\Big) N^{-\alpha} 
\sum\Sb a, r\endSb f(a)^2 a^{\alpha} \frac{f(r)r^{\alpha}}{\sqrt{r}} 
\\
&\ll (\log X) N^{-\alpha} \prod_p \Big( 1+\frac{f(p)}{\sqrt{p}}\Big)\Big( 
1+f(p)^2 p^{\alpha} + \frac{f(p)p^{\alpha}}{\sqrt{p}}\Big).
\\
\endalign
$$
The ratio of the above quantity to that in (15) is 
$$
\ll \exp\Big( -\alpha \log N + \sum_{p} \Big(f(p)^2 (p^{\alpha}-1) + 4f(p)p^{\alpha-\frac 12} \Big)\Big)
\ll \exp\Big(- \alpha\frac{\log N}{\log \log N}\Big).
$$
Combining this estimate with (14) and (15) we conclude that 
$$
M_2(R,X) \sim C_1 X (\log X) \prod_{p} \Big( 1+f(p)^2 -2\frac{f(p)}{\sqrt{p}}\Big).  
$$
This completes the proof of the Proposition.
\enddemo
\head 4. Extreme values of $L$-functions of cusp forms: Proof of Theorem 3  \endhead
\def\lam{\lambda}

\noindent Given $f \in H_k$ we define
$$
\omega(f):= \frac 3{\pi} \frac{(4\pi)^k}{\Gamma(k)} || f||^2, \tag{16a}
$$
where $||f||^2 = <\!f,f\!> = \int_{\Gamma\backslash {\Bbb H}}
 y^k |f(z)|^2 \frac{dx dy}{y^2}$ is the Petersson norm of $f$.
The weights $\omega(f)$ are related to the value at $1$ of 
the symmetric square $L$-function of $f$.  Namely, (see Iwaniec [7] for 
example)
$$
\omega(f) = L(1,\text{sym}^2 f)/\zeta(2). \tag{16b}
$$
We also know that the weights $\omega(f)$ are roughly of constant size; precisely,
$$
(\log k)^{-2} \ll \omega(f) \ll (\log k)^2. \tag{16c}
$$
For any two integers $m, \ n\ge 1$ we have
$$
\frac{12}{k-1} \sum_{f \in H_k} \frac{\lam_f(m)\lam_f(n)}{\omega(f)}
= \delta_{m,n} + 2\pi i^k \sum_{c=1}^{\infty} \frac{S(m,n;c)}{c}
J_{k-1}\Big(\frac{4\pi\sqrt{mn}}{c}\Big), \tag{17}
$$
where $\delta_{m,n}=1$ or $0$ depending on whether $m=n$ or not,
$J_{k-1}$ is the usual Bessel function, and $S(m,n;c)
= \sum_{a\pmod c}^{*} e(\tfrac{am+\overline{a}n}{c})$ is
Kloosterman's sum.  This is Petersson's formula, see Iwaniec [7].

If $x\le 2k$ then 
$$
|J_{k-1}(x)| = 
\Big|\sum_{\ell =0}^{\infty} \frac{(-1)^{\ell}}{\ell! (\ell +k-1)!} 
\Big(\frac x2\Big)^{2\ell +k-1} \Big|
\le \frac{(x/2)^{k-1}}{(k-1)!} \sum_{\ell=0}^{\infty} 
\frac{(x/2)^{2\ell}}{\ell! k^\ell} \le \frac{e^{x/2}(x/2)^{k-1}}{(k-1)!}. 
$$
Using this together with $|S(m,n;c)|\le c$ we obtain 
that if $4\pi \sqrt{mn} \le k/10$ then
$$
\align
\frac{12}{k-1} \sum_{f \in H_k} \frac{\lam_f(m)\lam_f(n)}{\omega(f)}
&= \delta_{m,n} 
+ O\Big( \frac{ e^{2\pi \sqrt{mn}} }{(k-1)!} 
\sum_{c=1}^{\infty} \Big(\frac{2\pi \sqrt{mn}}{c}\Big)^{k-1} \Big) 
\\
&=\delta_{m,n} + O(e^{-k}). \tag{18}\\
\endalign
$$

Let $r(n)$ be arbitrary real numbers and 
consider the resonator $R(f)=\sum_{n\le N} \lambda_f(n) r(n)$.  
If $N\le k/(40 \pi)$ 
then we obtain from (18) that 
$$
\frac{12}{k-1}\sum_{f \in H_k} \frac{R(f)^2}{\omega(f)} = 
\sum_{m,n \le N} r(m)r(n) ( \delta_{m,n} +O(e^{-k})) 
= \sum_{n\le N} r(n)^2 (1+O(ke^{-k})), \tag{19}
$$
where the last equality follows from Cauchy's inequality.  

Next we want to calculate the weighted average of 
$|R(f)|^2 L(\frac 12,f)$.  To do this we require 
an ``approximate functional equation'' for $L(\tfrac 12,f)$ 
which we now describe briefly.  We consider, for some
$c>\frac 12$,  
$$
\frac{1}{2\pi i}\int_{c-i\infty}^{c+i\infty} 
(2\pi)^{-s} \frac{\Gamma(s+\tfrac k2)}
{\Gamma(\tfrac k2)} 
L(s+\tfrac 12,f) \frac{ds}{s}. \tag{20a}
$$
We move the line of integration to the line Re$(s) =-c$ and use the 
functional equation.  The pole at $s=0$ leaves the residue $L(\tfrac 12,f)$ 
and thus (20a) equals 
$$
L(\tfrac 12,f) +\frac{i^k}{2\pi i} \int_{-c-i\infty}^{-c+i\infty} 
(2\pi)^{-s} \frac{\Gamma(-s+\frac k2)}{\Gamma(\frac k2)} 
L(-s+\frac 12,f) \frac{ds}{s}. 
$$
Replacing $-s$ by $s$ we deduce that 
$$
L(\tfrac 12,f) = (1+i^k) \frac{1}{2\pi i}\int_{c-i\infty}^{c+i\infty} 
(2\pi)^{-s} \frac{\Gamma(s+\tfrac k2)}
{\Gamma(\tfrac k2)} 
L(s+\tfrac 12,f) \frac{ds}{s}. \tag{20b}
$$
Defining, for real numbers $x> 0$, 
$$
V(x):= \frac{1}{2\pi i} \int_{c-i\infty}^{c+i\infty} 
(2\pi)^{-s} \frac{\Gamma(s+\frac k2)}{\Gamma(\frac k2)} 
x^{-s} \frac{ds}{s}, \tag{20c}
$$ 
and expanding $L(s+\frac 12,f)$ into its Dirichlet series we deduce 
from (20b) that 
$$
L(\tfrac 12,f) = (1+i^k) \sum_{n=1}^{\infty} 
\frac{\lam_f(n)}{\sqrt{n}} V(n).
\tag{20d}
$$
Moving the line of integration in (20c) to $c=k/2$ and $c=1-k/2$ 
we obtain respectively that
$$
V(x) \ll \Big(\frac{k}{2\pi x}\Big)^{\frac{k}{2}}, \qquad 
\text{and }\qquad V(x) = 1+ O\Big(\frac{(2\pi x)^{\frac k2 -1}}
{\Gamma(\frac k2)}
\Big). \tag{20e}
$$

Suppose now that $N\le \sqrt{k}/100$.  Then, using 
the Hecke relations and (18),  we obtain that 
$$
\align
\frac{12}{k-1} \sum_{f\in H_k} &\frac{R(f)^2}{\omega(f)} 
\sum_{r\le 2k} \frac{\lam_f(r)}{\sqrt{r}} V(r) 
\\
&=\frac{12}{k-1} \sum_{f\in H_k}\frac{1}{\omega(f)} 
\sum_{m,n \le N} r(m) r(n) \sum_{d|(m,n) }
\lam_f\Big(\frac{mn}{d^2}\Big) 
 \sum_{r\le 2k} \frac{\lam_f(r)}{\sqrt{r}} V(r)
\\
&=\sum_{m,n \le N} r(m) r(n)\sum_{d|(m,n) } \Big(\frac{d}{\sqrt{mn}}
V\Big(\frac{mn}{d^2}\Big)
+O(ke^{-k})\Big)\\
&= \sum_{m,n \le N} r(m) r(n)\frac{\sigma((m,n))}{\sqrt{mn}} 
+ O\Big(k^3 e^{-k} \sum_{n\le N} r(n)^2\Big).\\
\endalign
$$
The final inequality above follows upon using (20e) to replace
$V(mn/d^2)$ by $1$, and then noting that $\sum_{d|(m,n)} 1 \le k$ 
and that $\sum_{m,n\le N} |r(m)r(n)| \le M \sum_{n\le N}r(n)^2$ 
by Cauchy's inequality.  
Now suppose that $k\equiv 0 \pmod 4$.  Note that by (20e) 
the terms $n>2k$ contribute an amount $O(e^{-k})$ to $L(\tfrac 12,f)$.  
Therefore we deduce that 
$$
\frac{12}{k-1} \sum_{f\in H_k} \frac{R(f)^2}{\omega(f)} L(\tfrac 12,f) 
= \sum_{m,n \le N} r(m) r(n)\frac{\sigma((m,n))}{\sqrt{mn}} 
+ O\Big(k^3 e^{-k} \sum_{n\le N} r(n)^2\Big). \tag{21}
$$

To produce large values of $L(\tfrac 12,f)$ we choose $N=\sqrt{k}/100$, 
and choose the resonator coefficients $r(n)$ to be $f(n)$, where $f$ is 
the multiplicative function used in the proof of the lower bound in Theorem 2.1.  
Using Rankin's trick in (19) we obtain that 
$$
\frac{12}{k-1} \sum_{f\in H_k} \frac{R(f)^2}{\omega(f)} 
\sim \prod_{p} (1+f(p)^2). 
$$
Further, Rankin's trick and (21) give 
$$
\frac{12}{k-1} \sum_{f\in H_k} \frac{R(f)^2}{\omega(f)} L(\tfrac 12,f) 
\sim \prod_p \Big( 1+ f(p)^2  \Big( 1+\frac 1p\Big)+ 2\frac{f(p)}{\sqrt{p}}\Big),
$$
and the conclusion of Theorem 3 regarding large values follows.  
To obtain the conclusion concerning small values, we choose 
$r(n)$ to be $\mu(n)f(n)$.


\Refs
 \frenchspacing
 \widestnumber\key{20}
 


\ref\key 1
\by R. Balasubramanian 
\paper On the frequency of Titchmarsh's phenomenon for 
$\zeta(s)$-IV 
\jour Hardy-Ramanujan J.  
\vol 9 
\yr 1986 
\pages 1--10 
\endref

\ref\key 2 
\by R. Balasubramanian and K. Ramachandra 
\paper  On the frequency of Titchmarsh's phenomenon for 
$\zeta(s)$-III
\jour Proc. Indian Acad. Sci. 
\vol 86
\yr 1977
\pages 341--351
\endref 

\ref\key 3 
\by J.B. Conrey and K. Soundararajan 
\paper Real zeros of quadratic Dirichlet $L$-functions
\jour Invent. Math. 
\vol 150
\yr 2002
\pages 1--44
\endref

\ref\key 4
\by D.W. Farmer, S.M. Gonek and C.P. Hughes
\paper The maximum size of $L$-functions
\jour J. Reine Angew. Math. 
\vol 609
 \yr 2007
 \pages 215--236
 \endref

\ref\key  5
\by A. Granville and K. Soundararajan
\paper Large character sums
\jour J. Amer. Math. Soc. 
\vol 14
\pages 365--397
\yr 2001
\endref

\ref\key 6 
\by J. Hoffstein and P. Lockhart
\paper Omega results for automorphic $L$-functions
\inbook Automorphic forms, automorphic representations, and arithmetic
\publ AMS Proc. Symp. Pure Math.  66, part 2
\pages 239--250
\yr 1999
\endref

\ref\key 7 
\by H. Iwaniec
\book Topics in classical automorphic forms
\publ AMS Graduate Studies in Math. 17
\yr 1997
\pages xii+ 259
\endref


\ref\key 8 
\by D. Milicevic
\paper Large values of eigenfunctions on arithmetic hyperbolic manifolds
\jour Ph. D. thesis, Princeton University 
\yr 2006
\endref

\ref \key 9
\by H.L. Montgomery 
\paper Extreme values of the Riemann zeta function 
\jour Comment. Math. Helv. 
\vol 52 
\yr 1977
\pages 511--518 
\endref

\ref\key 10
\by N. Ng 
\paper Extreme values of $\zeta^{\prime}(\rho)$
\pages 16 pp.,
 e-print available at {\tt http://arxiv.org/abs/0706.1765}
\endref

\ref\key 11
\by Z. Rudnick and K. Soundararajan 
\paper Lower bounds for moments of $L$-functions
\jour Proc. Natl. Acad. Sci. USA 
\vol 102
\yr 2005
\pages 6837--6838
\endref

\ref\key 12
\by Z. Rudnick and K. Soundararajan 
\paper Lower bounds for moments of $L$-functions: symplectic 
and orthogonal examples
\inbook Multiple Dirichlet series, automorphic forms, and analytic number 
theory 
\publ AMS Proc. Sympos. Pure Math., 75
\yr 2006 
\pages 293--303
\endref

\ref \key 13
\by K. Soundararajan 
\paper Non-vanishing of quadratic Dirichlet $L$-functions 
at $s=1/2$ 
\jour Annals of Math. 
\vol 152 
\yr 2000 
\pages 447-488
\endref 

\ref\key 14
\by K. Soundararajan 
\paper Moments of the Riemann zeta-function 
\pages 11 pp
\jour Annals of Math.
 to appear,   e-print available at {\tt http://arxiv.org/abs/math/0612106}
\endref

\ref\key 15
\by E.C. Titchmarsh
\book The theory of the Riemann zeta-function 
\bookinfo Second Edition
\publ Oxford Univ. Press
\yr 1986 
\publaddr New York
\endref 

\endRefs
\enddocument